# A game theoretic model of preventive demographic measures

O.A. Malafeyev,

Saint-Petersburg State University

malafeyevoa@mail.ru

M.A. Suvorov,

Saint-Petersburg State University

st051952@student.spbu.ru

**Abstract**

The article constructs a game-theoretic model of preventive measures for the antagonistic game with nature, that is such a model in which two participants take part, one of which is a regulating body (government) that makes legislative decisions depending on the emerging external and internal conditions affecting the course of the demographic process, the second participant – internal and external conditions of the environment (nature). An example of such a game is given and the optimal strategy for the government is calculated.

## 1 Introduction

Virtually all countries at different periods of their history have faced the phenomenon of demographic crisis, and the Russian Federation was not left out of this trend. The last period of marked decline in the birth rate in Russia was in the 1990s and 2000s, which was associated with major changes in the post-Soviet space.

There is also a gradual process of change in the population structure from a traditional society with high birth and death rates to an industrialised and eventually post-industrial society with low birth and death rates.

The demographic situation is influenced by various factors, including the level of education, quality of life and socio-economic environment. In Russia, as in many other countries, demographic support programs are being introduced to regulate fertility. It is important to note that the impact of these programs on fertility can be diverse and depends on various factors.

However, all of them are aimed at stimulating the population to take such an important and responsible step as having children.

Let us consider what kind of demographic incentive programs can be disseminated among the population to improve fertility:

1. The program of maternal (family) capital. It is oriented towards providing money to families in the context of the birth or adoption of a first or second child. These funds, often paid in two instalments over a period of time, provide an opportunity for families to improve housing conditions or pay the education.

2. Mortgage loan packages for young families provide an opportunity to obtain favourablemortgage loan terms for the purchase of primary housing. This provides an opportunity to reduce the financial burden on young families, contributing to their stability during important life events.

3. Benefits and allowances for large families include a wide range of government incentives suchas transport allowances, free school meals, utility subsidies and other forms of financial support aimed at promoting the well-being of large families.

4. The development of preschool education is supported by the State in the form of the construction of new preschool establishments and the improvement of the qualifications of teaching staff. These measures are aimed at easing the responsibility of families and ensuring high quality education for young children.

5. State medical support covers programs to reduce maternal and child mortality rates byproviding free medical services for pregnant women and systematic newborn screening. These initiatives contribute to improving the overall health of mothers and children, which is strategically important for the well-being of families and society as a whole.

## 2 Antagonistic Problem Formulation

Suppose that the government needs to determine how population support programs will affect negative demographic trends. The uncertainty of long-run forecasts of external conditions makes a game-theoretic approach useful if we take the government and environmental conditions as agents I and II, respectively, who pursue directly opposite goals. Then, for a finite set of alternatives available to the government and a finite set of forecasts of environmental conditions, we can construct positional games.

Consider when one agent (the government) can undertake one of $n$ events, with event $j$ reducing the losses from phenomenon $j$ to zero and affecting the losses caused by action $i$ ($i \neq j$). Let the outcome from events with flow $j$ be $q_j > 0$, and the losses from phenomenon $i$, if a change $j$ ($j \neq i$) is made, be $r_{ij} \geqslant 0$. If the utility of nature (agent I) is understood as the loss it causes to the party in question (agent II), then this conflict can be modelled by a finite antagonistic game defined by a matrix of payoffs

$$H = \begin{pmatrix} q_1 & r_{12} + q_2 & \cdots & r_{1n} + q_n \\ r_{21} + q_1 & q_2 & \cdots & r_{2n} + q_n \\ \cdots & \cdots & \cdots & \cdots \\ r_{n1} + q_1 & r_{n2} + q_2 & \cdots & q_n \end{pmatrix}$$

Obviously, agent II will try to minimise his losses and, counting on the worst case, must assume that agent I (nature) is antagonistic to him. Let us consider some cases of antagonistic games.

# 3 Solution of the task at hand

## 3.1 Neutralising strategies of agents

Consider the case of independence. Here we assume that only event $j$ neutralises nature's strategy $j$, i.e. $r_{jj} = 0$, and in the other cases the loss is independent of whether the event is produced or not. For nature's strategy $j$, we will put the loss as in the previous case equal to $t_j$. For convenience, we will assume that $t_1 > t_2 > ... > t_n$. In addition, we will put $q_i = 0$. This case models conflicts in which the cost of additional interventions is negligible compared to the possible damage. Then the matrix $H$ takes the form

$$\begin{pmatrix} 0 & t_1 & t_1 & ... & t_1 \\ t_2 & 0 & t_2 & ... & t_2 \\ . & . & . & . & . \\ t_n & t_n & t_n & ... & 0 \end{pmatrix}$$

The game with the winning matrix $H$ of the previous case is reduced to the same form on the basis of the theorem on mixed strategies of players if all $p$ are equal to each other. Let us find the solution of the game under consideration. Let the following inequalities hold for all

$3 \leqslant k \leqslant n$

$$t_k > (k-2)/\sum_{i=1}^{k=1} \frac{1}{t_i}. \tag{1.1}$$

Then for all $j = 1, 2, ..., n$ the following inequalities hold

$$t_j > (n-1)/\sum_{i=1}^{n} \frac{1}{t_i} \tag{1.2}$$

as the numbers $l_j$ decrease with increasing index. Due to inequalities (1.2), the strategy $Y^*$, whose components are determined by the following equations

$$\eta_j^* = \frac{t_j \sum_{i=1}^{n} 1/t_i - (n-1)}{t_j \sum_{i=1}^{n} 1/t_i} \tag{1.3}$$

is the mixed strategy of agent II. Applying it against agent I's pure strategy $i$, he gets the value of

$$H(i, Y^*) = i_i \sum_{j \neq i} \eta_j^* = \left(1 - \eta_i^i\right) t_i = (n-1)/\sum_{k=1}^{n} \frac{1}{t_k}; \quad 1 \leqslant i \leqslant n. \tag{1.4}$$

If agent I applies the strategy $X^*$ for which

$$\xi_i^* = \left(t_i \sum_{k=1}^{n} \frac{1}{t_k}\right)^{-1} \tag{1.5}$$

against the pure strategy $j$ of agent II, then he will get

$$H(X_2^* j) = \sum_{i \neq j} t_i \xi_i^* = \frac{n-1}{\sum_{k=1}^{n} 1/t_k}, \quad 1 \leqslant j \leqslant n. \tag{1.6}$$

It follows from (1.4) and (1.6) that the value $(n-1)\left(\sum_{k=1}^{n} 1/t_k\right)^{-1}$ is the value of the game, and the mixed strategies $X^*$ and $Y^*$ are the optimal mixed strategies of agents I and II.

Consider an example where the support measures add $t_1 = 30$ thousand people, $t_2 = 28$ thousand people, $t_3 = 26$ thousand people, and $t_4 = 24$ thousand people, $t_5 = 22$ thousand people, respectively. Taking the state support measures for agent II and evaluating the utility of agent I (nature) by the loss $t_i$, we arrive at an antagonistic finite game, which models the considered conflict with nature and is given by the gain matrix

$$H = \begin{pmatrix} 0 & 30 & 30 & 30 & 30 \\ 28 & 0 & 28 & 28 & 28 \\ 26 & 26 & 0 & 26 & 26 \\ 24 & 24 & 24 & 0 & 24 \\ 22 & 22 & 22 & 22 & 0 \end{pmatrix}$$

It is easy to establish that this game belongs to the case of independence of the above example, and therefore for its solution we check the fulfilment of inequalities (1.2). These inequalities are fulfilled because inequalities (1.1) are true. Consequently, we can use formula (1.3) to compute the optimal strategy of the state.

As a result, we find $Y^* = (0.315, 0.266, 0.21, 0.143, 0.066)$, and calculate the value of the game by the formula

$$v = \frac{n-1}{\sum_{k=1}^{n} \frac{1}{t_k}} = 20.8$$

Thus, the optimal strategy of the government is to allocate a share $y_i$ to the $i$ program, and the value of the game (i.e. the expected income in the least favourable case) will be $v = 20.8*10^3$ newborns.

## 3.2  Mixed strategies of agents

Suppose now that inequalities (1.1) are not fulfilled for some values of $k$. Then there exists a minimum index $\omega$ for which the inequalities (1.1) are still satisfied, and for $\omega + 1$ – are $i_k > (k-2) \left( \sum_{i=1}^{k-1} \frac{1}{t_i} \right)^{-1}, 3 \leqslant k \leqslant \omega_2$ violated, i.e.

(1.7)

$$t_{\omega+1} \leqslant (\omega - 1) \left( \sum_{i=1}^{\omega} \frac{1}{t_i} \right)^{-1}$$

(1.8)

By virtue of inequalities (1.7) the numbers

$$\eta_j^* = \frac{t_i \sum_{i=1}^{\omega} 1/t_i - (\omega - 1)}{t_j \sum_{i=1}^{\omega} 1/i_j}, \quad 1 \leqslant j \leqslant \omega,$$

(1.9)

are non-negative and sum to one. Then, the strategy $Y^* = (\eta_1*, ..., \eta_0*, 0, ..., 0)$ is determined by (1.9) and is a mixed strategy of agent II. Applying it against agent I's pure strategy $i$, he will lose $H(i, Y^*)$, and for $i = 1, 2, ..., \omega$ the following is true

$$H(i, Y^*) = (1 - \eta_i^*) t_i = (\omega - 1) \left( \sum_{k=1}^{\omega} \frac{1}{t_k} \right)^{-1}, \tag{1.10}$$

and for $i \geqslant \omega + 1$ from the fact that the numbers $t_i$ are decreasing. It follows from inequalities (1.8) that.

$$H(i, Y^*) = (1 - \eta_i^*) t_i = t_i \leqslant t_{\omega+1} \leqslant (\omega - 1) \left( \sum_{k=1}^{\omega} \frac{1}{t_k} \right)^{-1}, \tag{1.11}$$

Applying the strategy $X^* = (\xi_1^*, \xi_2^*, \ldots, \xi_\omega^*, 0, \ldots, 0)$, where

$$\xi_i^* = \begin{cases} \left( t_i \sum_{k=1}^{j} \frac{1}{t_k} \right)^{-1}, & \text{for } 1 \leqslant i \leqslant \omega \\ 0, & \text{for } i > \omega \end{cases}$$

(1.12)

and for $1 \leqslant j \leqslant \omega$, agent l gets

Thus, it follows from (1.10) to (1.14) that

$$H(i, Y^*) \leqslant H(X^*, Y^*) \leqslant H(X^*, j), \tag{1.15}$$

$$1 \leqslant i \leqslant n, \quad 1 \leqslant j \leqslant n,$$

where

$$H(X^*, j) = \sum_{i \neq j}^{\omega} t_i \xi_i^* = (\omega - 1) \left( \sum_{k=1}^{\omega} \frac{1}{t_k} \right)^{-1} \tag{1.13}$$

If $j > \omega$, then

$$H(X^*, j) = \sum_{i=1}^{\omega} t_i \xi_i^* = \omega \left( \sum_{k=1}^{\omega_j} \frac{1}{t_k^*} \right)^{-1} > (\omega - 1) \left( \sum_{k=1}^{\omega} \frac{1}{t_k} \right)^{-1} \tag{1.14}$$

$$H(X^*, Y^*) = (\omega - 1) \left( \sum_{k=1}^{\omega} \frac{1}{t_k} \right)^{-1} \tag{1.16}$$

Consequently,

$$v = (\omega - 1) \left( \sum_{k=1}^{\omega} \frac{1}{t_k} \right)^{-1}, \tag{1.17}$$

and optimal strategies of agents are defined by the equations (1.9) and (1.12).

Suppose, the losses $t_i$ (i = 1, 2, ..., 5) take such values that the payoff matrix of agent I has the form

$$H_1 = \begin{pmatrix} 0 & 1/3 & 1/3 & 1/3 & 1/3 \\ 1/35 & 0 & 1/35 & 1/35 & 1/35 \\ 1/4 & 1/4 & 0 & 1/4 & 1/4 \\ 1/45 & 1/45 & 1/45 & 0 & 1/45 \\ 1/5 & 1/5 & 1/5 & 1/5 & 0 \end{pmatrix}$$

Checking inequalities (1.1) shows that the minimum of the indices of $\omega$ for which these inequalities are still satisfied is 3, while for $\omega + 1 = 4$ they are violated. Consequently, to calculate the components of the vector of the optimal strategy of the state (agent II) we need to use formula (1.3), and the value of the game should be found by formula (1.16). Having performed the calculations, we obtain $Y^* = (0.55, 0.31, 0.14, 0, 0, 0)$. And the value of the game will be $v = 0.1656$.

## Conclusions

As a result of this work, two problems of the antagonistic game were considered - cases of neutralising and mixed strategy of agents. With the help of modern computing machines, the calculation was made and the following conclusions were obtained: the optimal strategy of the government is to under $i$ program, and the value of the game (i.e. the expected income in the least favourable case) will be $v = 20.8 * 10^3$ newborns. In this case, optimal resource allocation leads to an improvement in the situation, providing ample prospects for increasing fertility and hence improving the demographic situation of the country. The obtained results allow the government to make informed decisions, taking into account the least favourable variant.

## References

1. Malafeyev O.A., Nemnyugin S.A. Stohasticheskaya model' social'no-ekonomicheskoj dinamiki // V sbornike: Ustojchivost' i processy upravleniya. Materialy III mezhdunarodnoj konferencii - 2015 - S. 433-434.

2. Zaitseva I., Poddubnaya N., Malafeyev O., Vanina A., Novikova E. Solving a dynamic assignment problem in the socioeconomic system // Journal of Physics: Conference Series – 2019 – S. 012092.

3. Malafeyev O., Zaitseva I., Pavlov I., Shulga A., Sychev S. Badin G. company life cycle model: the influence of interior and exterior factors // AIP Conference Proceedings – 2020 - S. 420027.

4. Malafeyev O.A., Rylow D., Pichugin Y.A., Zaitseva I. A statistical method for corrupt agents detection // AIP Conference Proceedings – 2018 - S. 100014.

5. Malafeyev O.A., Rylow D., Zaitseva I., Ermakova A., Shlaev D., Multistage voting model with alternative elimination // AIP Conference Proceedings – 2018 - S. 100012.

6. Malafeyev O.A., Redinskikh N.D., Nemnyugin S.A., Kolesin I.D., Zaitseva I.V. The optimization problem of preventive equipment repair planning // AIP Conference Proceedings – 2018 S. 100013.

7. Kolesin I., Malafeyev O., Andreeva M., Ivanukovich G. Corruption: taking into account the psychological mimicry of officials // AIP Conference Proceedings – 2017 - S. 170014.

8. Malafeyev O., Rylow D., Novozhilova L., Zaitseva I., Popova M., Zelenkovskii P. Gametheoretic model of dispersed material drying process // AIP Conference Proceedings - 2017 - S. 020063.

9. Malafeyev O., Lakhina J., Redinskikh N., Smirnova T., Smirnov N., Zaitseva I. A mathematical model of production facilities location // Journal of Physics: Conference Series – 2019 - S. 012090.

10. Zaitseva I., Ermakova A., Shlaev D., Malafeyev O., Strekopytov S. Game-theoretical model of labour force training // Journal of Theoretical and Applied Information Technology – 2018 - Vol. 96 № 4 - S. 978-983.

11. Malafeyev O.A., Ahmadyshina A.R., Demidova D.A., Model' tendera na rynke rielterskih uslug s uchetom korrupcii // V knige: Stroitel'stvo i ekspluataciya energoeffektivnyh zdanij (teoriya i praktika s uchetom korrupcionnogo faktora) (Passivehouse) – 2015 - S. 161-168.

12. Malafeyev O.A., Redinskih N.D., Stohasticheskoe ocenivanie i prognoz effektivnosti razvitiya firmy v usloviyah korrupcionnogo vozdejstviya // V sbornike: Ustojchivost' i process upravleniya. Materialy III mezhdunarodnoj konferencii – 2015 - S. 437-438.

13. Titarenko M.L., Ivashov L.G., Kefeli I.F., Malafeyev O.A., i dr. Evrazijskaya duga nestabil'nosti i problemy regional'noj bezopasnosti ot Vostochnoj Azii do Severnoj Afriki, Kollektivnaya monografiya - SPb: Studiya NP-Print – 2013 - 576 s.

14. Malafeyev O.A. The existence of situations of $\epsilon$-equilibrium in dynamic games with dependent movements// USSR Computational Mathematics and Mathematical Physics – 1974 - Vol. 14. № 1 - S. 88-99.

15. Malafeyev O.A., Kolesin I.D., Kolokoltsov V.N.. i dr. Vvedenie v modelirovanie korrupcionnyh sistem i processov, kollektivnaya monografiya. Tom 1. - Stavropol': Izdatel'skij dom "Tesera" – 2016 - 224 s.

16. Malafeyev O., Kupinskaya A., Awasthi A., Kambekar K.S. Random walks and market efficiency in chinese and indian markets // Statistics, Optimization and Information Computing – 2019 - T. 7. № 1. - S. 1-25.

17. Marahov V.G., Malafeyev O.A. Dialog filosofa i matematika: «O filosofskih aspektah matematicheskogo modelirovaniya social'nyh preobrazovanij XXI veka» // V sbornike: Filosofiya poznaniya i tvorchestvo zhizni. Sbornik statej – 2014 - S. 279292.

18. Zaitseva I., Dolgopolova A., Zhukova V., Malafeyev O., Vorokhobina Y. Numerical method for computing equilibria in economic system models with labor force // AIP Conference Proceedings – 2019 - S. 450060.

19. Malafeyev O.A. On the existence of Nash equilibria in a noncooperative n-person game with measures as coefficients // Communications in Applied Mathematics and Computational Science – 1995 - Vol. 5. № 4 - S. 689-701.

20. Malafeyev O.A., Redinskih N.D., Smirnova T.E. Setevaya model' investirovaniya proektov s korrupciej // Processy upravleniya i ustojchivost' – 2015 - Vol. 2. № 1 - S. 659-664.

21. Kirjanen A.I., Malafeyev O.A., Redinskikh N.D. Developing industries in cooperative interaction: equilibrium and stability in process with lag // Statistics, Optimization and Information Computing – 2017 – Vol. 5. № 4 - S. 341-347.

22. Troeva M.S., Malafeyev O.A. Ravnovesie v diffuzionnoj konfliktnoj modeli ekonomiki so mnogimi uchastnikami // V sbornike: Dinamika, optimizaciya, upravlenie. Ser. "Voprosy mekhaniki i processov upravleniya" – 2004 - S. 146-153.

23. Pichugin Yu.A., Malafeyev O.A. Ob ocenke riska bankrotstva firmy// V knige: Dinamicheskie sistemy: ustojchivost', upravlenie, optimizaciya. Tezisy dokladov – 2013 - Vol. 204-206.

24. Malafeyev O.A., Redinskih N.D., Alferov G.V., Smirnova T.E. Korrupciya v modelyah aukciona pervoj ceny // V sbornike: Upravlenie v morskih i aerokosmicheskih sistemah (UMAS2014) – 2014 - S. 141-146.

25. Kefeli I.F., Malafeyev O.A. O matematicheskih modelyah global'nyh geopoliticheskih processov mnogoagentnogo vzaimodejstviya // Geopolitika i bezopasnost' – 2013 - № 2 (22) - S. 44-57.

26. Kolokoltsov V.N., Malafeyev O.A. Corruption and botnent defence: a mean field game approach // International Journal of Game Theory – 2018 - Vol. 47. № 3 - S. 977-999.

27. Malafeyev O.A., Redinskikh N.D. Compromise solution in the problem of change state control for the material body exposed to the external medium // AIP Conference Proceedings – 2018 - S. 080017.

28. Malafeyev O.A., Rylow D., Kolpak E.P., Nemnyugin S.A., Awasthi A. Corruption dynamics model // AIP Conference Proceedings – 2017 - S. 170013.

29. Zaitseva I., Ermakova A., Shlaev D., Malafeyev O., Kolesin I., Modeling of the labour force redistribution in investment projects with account of their delay // IEEE International Conference on Power, Control, Signals and Instrumentation Engineering – 2017 - S. 68-70.

30. Malafeyev O., Onishenko V., Zubov A., Bondarenko L., Orlov V., Petrova V., Kirjanen A., Zaitseva I. Optimal location problem in the transportation network as an investment project: a numerical method // AIP Conference Proceedings – 2019 - S. 450058.

31. Vlasov M.A., Glebov V.V., Malafeyev O.A., Novichkov D.N. Experimental study of an electron beam in drift space // Soviet Journal of Communications Technology and Electronics – 1986 Vol. 31. №3 - S. 145.

32. Malafeyev O.A. O dinamicheskih igrah s zavisimymi dvizheniyami // Doklady Akademii nauk SSSR – 1973 - Vol. 213. № 4 - S. 783-786.

33. Malafeyev O.A., Redinskih N.D., Alferov G.V., Model' aukciona s korrupcionnoj komponentoj // Vestnik Permskogo universiteta. Seriya: Matematika. Mekhanika. Informatika – 2015 - № 1 (28) - S. 30-34.

34. Malafeyev O.A., Konfliktno upravlyaemye processy so mnogimi uchastnikami, avtoreferat dis. doktora fizikomatematicheskih nauk / LGU im. A. A. ZHdanova. Leningrad – 1987 - 44 s.

35. Ivanyukovich G.A., Malafeyev O.A., Zaitseva I.V., Kovshov A.M., Zakharov V.V., Zakharova N.I. To the evaluation of the parameters of the regression equation between the radiometric and geological testing // JOP Conference Series: Metrological Support of Innovative Technologies - 2020 - S. 32079.

36. Malafeyev O.A., Sushchestvovanie situacii ravnovesiya v beskoalicionnyh differencial'nyh igrah dvuh lic s nezavisimymi dvizheniyami // Vestnik Leningradskogo universiteta. Seriya 1: Matematika, mekhanika, astronomiya – 1980 - № 4 - S. 1216.

37. Malafeyev O.A., Stohasticheskoe ravnovesie v obshchej modeli konkurentnoj dinamiki// V sbornike: Matematicheskoe modelirovanie i prognoz social'no-ekonomicheskoj dinamiki v usloviyah konkurencii i neopredelennosti - 2004-S. 143-154.

38. Malafeyev O.A., Eremin D.S. Konkurentnaya linejnaya model' ekonomiki // V sbornike: Processy upravleniya i ustojchivost' – 2008 - S. 425-435.

39. Malafeyev O.A., Redinskikh N.D. Quality estimation of the geopolitical actor development strategy// CNSA 2017 Proceedings – 2017- http://dx.doi.org/10.1109/CNSA.2017.7973986

40. Ivashov L.G., Kefeli I.F., Malafeyev O.A. Global'naya arkticheskaya igra i ee uchastniki// Geopolitika i bezopasnost' - 2014 - № 1 (25) - S. 34-49.

41. Kolokoltsov V.N., Malafeyev O.A., Mean field game model of corruption // arxiv.org
https://arxiv.org/abs/1507.03240

42. Zaitseva I.V., Malafeyev O.A., Zakharov V.V., Zakharova N.I., Orlova A.Yu. Dynamic distribution of labour resources by region of investment // Journal of Physics: Conference Series. – 2020 - S. 012073.

43. Malafeyev O.A. Sushchestvovanie situacij ravnovesiya v differencial'nyh beskoalicionnyh igrah so mnogimi uchastnikami // Vestnik Leningradskogo universiteta. Seriya 1: Matematika, mekhanika, astronomiya. – 1982 - № 13 - S. 4046.
10

44. Malafeyev O.A., Petrosyan L.A. Igra prostogo presledovaniya na ploskosti s prepyatstviem // Upravlyaemye sistemy – 1971 - № 9 - S. 31-42.

45. Malafeyev O. A, Galtsov M., Zaitseva I., Sakhnyuk P., Zakharov V., Kron R. Analysis of trading algorithms on the platform QIUK // Proceedings - 2020 2nd International Conference on Control Systems, Mathematical Modeling, Automation and Energy Efficiency – 2020 - S. 305-311.

46. Malafeyev O.A. Ustojchivye beskoalicionnye igry N lic // Vestnik Leningradskogo universiteta. Seriya 1: Matematika, mekhanika, astronomiya. – 1978 - № 4 - S. 55-58.

47. Malafeyev O.A., CHernyh K.S. Prognosticheskaya model' zhiznennogo cikla firmy v konkurentnoj srede // V sbornike: Matematicheskoe modelirovanie i prognoz social'noekonomicheskoj dinamiki v usloviyah konkurencii i neopredelennosti. – 2004 - S. 239-255.

48. Drozdov G.D., Malafeyev O.A. Modelirovanie tamozhennogo dela – SPb: Izd-vo SPbGUSE – 2013 - 255 s. 60. Pichugin YU.A., Malafeyev O.A., Alferov G.V. Ocenivanie parametrov v zadachah konstruirovaniya mekhanizmov robotov-manipulyatorov // V sbornike: Ustojchivost' i processy upravleniya. Materialy III mezhdunarodnoj konferencii – 2015 - S. 141-142.

49. Malafeyev O.A., Zajceva I.V., Komarov A.A., SHvedkova T.YU., Model' korrupcionnogo vzaimodejstviya mezhdu kommercheskoj organizaciej i otdelom po bor'be s korrupciej // V knige: Linejnaya algebra s prilozheniyami k modelirovaniyu korrupcionnyh sistem i processov – 2016 - S. 342-351.

50. Malafeyev O.A., Strekopytova O.S. Teoretiko-igrovaya model' vzaimodejstviya korrumpirovannogo chinovnika s klientom: psihologicheskie aspekty // V knige: Vvedenie v modelirovanie korrupcionnyh sistem i processov. Malafeyev O.A., i dr. Kollektivnaya monografiya. Pod obshchej redakciej d.f.-m.n., professora O. A. Malafeyeva - 2016. S. 134- 151.

51. Zaitseva I.V., Malafeyev O.A., Zakharov V.V., Smirnova T.E., Novozhilova L.M. Mathematical model of network flow control // IOP Conference Series: Materials Science and Engineering – 2020 - S. 012036.

52. Malafeyev O.A., Petrosyan L.A. Differencial'nye mnogokriterial'nye igry so mnogimi uchastnikami // Vestnik Leningradskogo universiteta. Seriya 1: Matematika, mekhanika, astronomiya – 1989 - № 3 - S. 27-31.

53. Malafeyev O.A., Kefeli I.F. Nekotorye zadachi obespecheniya oboronnoj bezopasnosti// Geopolitika i bezopasnost' – 2013 - № 3 (23) - S. 84-92.

54. Malafeyev O.A., Redinskih N.D. Stohasticheskij analiz dinamiki korrupcionnyh gibridnyh setej // V sbornike: Ustojchivost' i kolebaniya nelinejnyh sistem upravleniya (konferenciya Pyatnickogo) – 2016 - S. 249-251.

55. Malafeyev O.A., Farvazov K.M. Statisticheskij analiz indikatorov korrupcionnoj deyatel'nosti v sisteme gosudarstvennyh i municipal'nyh zakupok // V knige: Vvedenie v modelirovanie korrupcionnyh sistem i processov – 2016 - S. 209-217.

56. Malafeyev O.A., Zajceva I.V., Zenovich O.S., Rumyancev N.N., Grigor'eva K.V., Ermakova A.N., Rezen'kov D.N., SHlaev D.V. Model' raspredeleniya resursov pri vzaimodejstvii korrumpirovannoj struktury s antikorrupcionnym agentstvom // V knige: Vvedenie v modelirovanie korrupcionnyh sistem i processov – 2016 - S. 102-106.

57. Malafeyev O., Parfenov A., Smirnova T., Zubov A., Bondarenko L., Ugegov N., Strekopytova M., Strekopytov S., Zaitseva I. Game-theoretical model of cooperation between producers // V sbornike: AIP Conference Proceedings – 2019 - S. 450059.

58. Neverova E.G., Malafeyev O.A., Alferov G.V. Nelinejnaya model' upravleniya antikorrupcionnymi meropriyatiyami // V sbornike: Ustojchivost' i processy upravleniya. Materialy III mezhdunarodnoj konferencii – 2015 - S. 445-446.

59. Malafeyev O.A., Marahov V.G. Evolyucionnyj mekhanizm dejstviya istochnikov i dvizhushchih sil grazhdanskogo obshchestva v sfere finansovoj i ekonomicheskoj komponenty XXI veka // V sbornike: K. Marks i budushchee filosofii Rossii. Busov S.V., Dudnik S.I. i dr. - 2016. S. 112-135.

60. Malafeyev O.A., Borodina T.S., Kvasnoj M.A., Novozhilova L.M., Smirnov I.A. Teoretikoigrovaya zadacha o vliyanii konkurencii v korrupcionnoj srede // V knige: Vvedenie v modelirovanie korrupcionnyh sistem i processov – 2016 - S. 8896.

61. Malafeyev O.A., Zajceva I.V., Koroleva O.A., Strekopytova O.S. Model' zaklyucheniya kontraktov s vozmozhno korrumpirovannym chinovnikom-principalom // V knige: Vvedenie v modelirovanie korrupcionnyh sistem i processov – 2016 - S. 114-124.

62. Malafeyev O.A., Sajfullina D.A. Mnogoagentnoe vzaimodejstvie v transportnoj zadache s korrupcionnoj komponentoj // V knige: Vvedenie v modelirovanie korrupcionnyh sistem i processov – 2016 - S. 218-224.

63. Malafeyev O.A., Neverova E.G., Smirnova T.E., Miroshnichenko A.N. Matematicheskaya model' processa vyyavleniya korrupcionnyh elementov v gosudarstvennoj sisteme upravleniya // V knige: Vvedenie v modelirovanie korrupcionnyh sistem i processov – 2016 - S. 152-179.

64. Malafeyev O.A., Koroleva O.A., Neverova E.G. Model' aukciona pervoj ceny s vozmozhnoj korrupciej // V knige: Vvedenie v modelirovanie korrupcionnyh sistem i processov – 2016 - S. 96-102.

65. Malafeyev O.A., Neverova E.G., Petrov A.N. Model' processa vyyavleniya korrupcionnyh epizodov posredstvom inspektirovaniya otdelom po bor'be s korrupciej // V knige: Vvedenie v modelirovanie korrupcionnyh sistem i processov – 2016 - S. 106-114.

66. Malafeyev O.A., Parfenov A.P. Dinamicheskaya model' mnogoagentnogo vzaimodejstviya mezhdu agentami korrupcionnoj seti // V knige: Vvedenie v modelirovanie korrupcionnyh sistem i processov – 2016 - S. 55-63.

67. Malafeyev O.A., Andreeva M.A., Gus'kova YU.YU., Mal'ceva A.S. Poisk podvizhnogo ob"ekta pri razlichnyh informacionnyh usloviyah // V knige: Vvedenie v modelirovanie korrupcionnyh sistem i processov – 2016 - S.63-88.

68. Zajceva I.V., Popova M.V., Malafeyev O.A. Postanovka zadachi optimal'nogo raspredeleniya trudovyh resursov po predpriyatiyam s uchetom izmenyayushchihsya uslovij // V knige: Innovacionnaya ekonomika i promyshlennaya politika regiona (EKOPROM-2016) – 2016 - S. 439-443.

69. Malafeyev O.A., Novozhilova L.M., Redinskih N.D., Rylov D.S., Gus'kova YU.YU. Teoretikoigrovaya model' raspredeleniya korrupcionnogo dohoda // V knige: Vvedenie v modelirovanie korrupcionnyh sistem i processov – 2016 - S. 125-134.

70. Malafeyev O.A., Salimov V.A., SHarlaj A.S. Algoritm ocenki bankom kreditosposobnosti klientov pri nalichii korrupcionnoj sostavlyayushchej // Vestnik Permskogo universiteta. Seriya: Matematika. Mekhanika. Informatika – 2015 - № 1 (28) - S. 35-38.

71. Kefeli I.F., Malafeyev O.A. Problemy ob"edineniya interesov gosudarstv EAES, SHOS i BRIKS v kontekste teorii kooperativnyh igr // Geopolitika i bezopasnost' – 2015 - № 3 (31) S. 33-41.

72. Kefeli I.F., Malafeyev O.A., Marahov V.G. i dr. Filosofskie strategii social'nyh preobrazovanij XXI veka – SPb: Izd-vo SPbGU – 2014 - 144 s.

73. Malafeyev O.A., Andreeva M.A., Alferov G.V., Teoretiko-igrovaya model' poiska i perekhvata v N-sektornom regione ploskosti // Processy upravleniya i ustojchivost' – 2015 - Vol. 2. № 1 S. 652-658.

74. Zajceva I.V., Malafeyev O.A. Issledovanie korrupcionnyh processov i sistem matematicheskimi metodami // Innovacionnye tekhnologii v mashinostroenii, obrazovanii i ekonomike – 2017 -

Vol. 3. № 1-1 (3) - S. 7-13.

75. Kolchedancev L.M., Legalov I.N., Bad'in G.M., Malafeyev O.A., Aleksandrov E.E., Gerchiu A.L., Vasil'ev YU.G. Stroitel'stvo i ekspluataciya energoeffektivnyh zdanij (teoriya i praktika s uchetom korrupcionnogo faktora) (Passivehouse)- Borovichi: NP "NTO strojindustrii SanktPeterburga" - 2015 - 170 S.

76. Malafeyev O.A., Novozhilova L.M., Kvasnoj M.A., Legalov I.N., Primenenie metodov setevogo analiza pri proizvodstve energoeffektivnyh zdanij s uchetom korrupcionnogo faktora // V knige: Stroitel'stvo i ekspluataciya energoeffektivnyh zdanij (teoriya i praktika s uchetom korrupcionnogo faktora) (Passivehouse) – 2015 - S. 146-161.

77. Malafeyev O., Awasthi A., Zaitseva I., Rezenkov D., Bogdanova S. A dynamic model of functioning of a bank // AIP Conference Proceedings. International Conference on Electrical, Electronics, Materials and Applied Science – 2018 - S. 020042.

78. Malafeyev O.A. Obzor literatury po modelirovaniyu korrupcionnyh sistem i processov, ch.I // V knige: Vvedenie v modelirovanie korrupcionnyh sistem i processov – 2016 - S. 9-17.

79. Asaul A.N., Lyulin P.B., Malafeyev O.A. Matematicheskoe modelirovanie vzaimodejstvij organizacii kak zhivoj sistemy // Vestnik Hmel'nickogo nacional'nogo universiteta Ekonomicheskie nauki – 2013 - № 6-2 (206) - S. 215-220.

80. Awasthi A., Malafeyev O.A. Is the indian dtock market efficient – a comprehensive study of Bombay stock exchange indices // arxiv.org - https://arxiv.org/abs/1510.03704

81. Malafeyev O.A., Demidova D.A. Modelirovanie processa vzaimodejstviya korrumpirovannogo predpriyatiya federal'nogo otdela po bor'be s korrupciej // V knige: Vvedenie v modelirovanie korrupcionnyh sistem i processov – 2016 - S. 140-152.

82. Malafeyev O.A., Chernyh K.S. Matematicheskoe modelirovanie razvitiya kompanii // Ekonomicheskoe vozrozhdenie Rossii – 2005 - № 2 - S. 23.

83. Kulakov F.M., Alferov G.V., Malafeyev O.A. Kinematicheskij analiz ispolnitel'noj sistemy manipulyacionnyh robotov // Problemy mekhaniki i upravleniya: Nelinejnye dinamicheskie sistemy – 2014 - № 46 - S. 31-38.

84. Kulakov F.M., Alferov G.V., Malafeyev O.A. Dinamicheskij analiz ispolnitel'noj sistemy manipulyacionnyh robotov // Problemy mekhaniki i upravleniya: Nelinejnye dinamicheskie sistemy – 2014 - № 46 - S. 39-46.

85. Zajceva I.V., Malafeyev O.A., Stepkin A.V., CHernousov M.V., Kosoblik E.V. Modelirovanie ciklichnosti razvitiya v sisteme ekonomik // Perspektivy nauki – 2020 - № 10 (133) - S. 173-176.

98. Bure V.M., Malafeyev O.A., Some game-theoretical models of conflict in finance // Nova Journal of Mathematics, Game Theory, and Algebra. –1996 - T. 6. № 1 - S. 7-14.

86. Malafeyev O.A., Redinskih N.D., Parfenov A.P., Smirnova T.E. Korrupciya v modelyah aukciona pervoj ceny // V sbornike: Instituty i mekhanizmy innovacionnogo razvitiya: mirovoj opty i rossijskaya praktika. – 2014 - S. 250-253.

87. Neverova E.G., Malafeyev O.A. A model of interaction between anticorruption authority and corruption groups // AIP Conference Proceedings – 2015 - S. 450012.

88. Irina Zaitseva, Oleg Malafeyev, Yuliya Marenchuk, Dmitry Kolesov and Svetlana Bogdanova Competitive Mechanism For The Distribution Of Labor Resources In The Transport Objective // Journal of Physics: Conference Series – 2018 - Volume 1172 - https://doi.org/10.1088/17426596/1172/1/012089 .

89. Malafeyev O., Awasthi A., Zaitseva I., Rezenkov D., Bogdanova S. A dynamic model of functioning of a bank // AIP Conference Proceedings – 2018 - S. 020042.

90. Kolokoltsov V.N., Malafeyev O.A., Understanding game theory: introduction to the analysis of many agent systems with competition and cooperation - World Scientific Publishing Company.: 2010 – 286 c.

91. Zaitseva I., Malafeyev O., Strekopytov S., Bondarenko G., Lovyannikov D. Mathematical model of regional economy development by the final result of labour resources // AIP Conference Proceedings – 2018 - S. 020011.

92. Malafeyev O.A., Nemnyugin S.A., Ivaniukovich G.A. Stohastic models of social-economic dynamics // In Proc. International Conference "Stability and Control Processes" in Memory of V.I. Zubov (SCP) – 2015 - S. 483-485.

93. Malafeyev O.A., Redinskih N.D., Gerchiu A.L. Optimizacionnaya model' razmeshcheniya korrupcionerov v seti // V knige: Stroitel'stvo i ekspluataciya energoeffektivnyh zdanij (teoriya i praktika s uchetom korrupcionnogo faktora) (Passivehouse). - Borovichi: 2015. - S. 128-140.

94. Malafeyev O.A., Koroleva O.A., Vasil'ev YU.G. Kompromissnoe reshenie v aukcione pervoj ceny s korrumpirovannym aukcionistom// V knige: Stroitel'stvo i ekspluataciya energoeffektivnyh zdanij (teoriya i praktika s uchetom korrupcionnogo faktora) (Passivehouse). – Borovichi: 2015. - S. 119-127.

95. Malafeyev O.A., Koroleva O.A., Neverova E.G. Model' aukciona pervoj ceny s vozmozhnoj korrupciej // Vvedenie v modelirovanie korrupcionnyh sistem i processov. – Stavropol': 2016. - S. 96-102.

96. Drozdov G.D., Malafeyev O.A., Modelirovanie mnogoagentnogo vzaimodejstviya processov strahovaniya, monografiya. - SPb: Sankt-Peterburgskij gos. un-t servisa i ekonomiki, 2010. 179 s.

97. Malafeyev O., Saifullina D., Ivaniukovich G., Marakhov V., Zaytseva I. The model of multiagent interaction in a transportation problem with a corruption component // AIP Conference Proceedings – 2017 - S. 170015.

98. Malafeyev O., Farvazov K., Zenovich O., Zaitseva I., Kostyukov K., Svechinskaya T. Geopolitical model of investment power station construction project implementation // AIP Conference Proceedings – 2018 - S. 020066.